\documentclass[12pt]{amsart}
\usepackage{amsmath,amssymb,fullpage}

\newtheorem{thm}{Theorem}[section]
\newtheorem{cor}[thm]{Corollary}
\newtheorem{lem}[thm]{Lemma}

\numberwithin{equation}{section}

\newcommand{\beq}{\begin{equation}}
\newcommand{\eeq}{\end{equation}}

\newtheorem{rem}[thm]{Remark}

\newcommand{\diag}{\operatorname{diag}}

\newcommand{\row}{\operatorname{row}}
\newcommand{\col}{\operatorname{col}}

\numberwithin{equation}{section}

\newcommand{\bbm}{\begin{bmatrix}}
\newcommand{\ebm}{\end{bmatrix}}

\title[Contractive determinantal
representations of stable polynomials]
{Matrix-valued Hermitian positivstellensatz, lurking contractions, and contractive determinantal
representations of stable polynomials}

\author[Grinshpan]{Anatolii Grinshpan}
\address{Department of Mathematics \\
Drexel University\\
3141 Chestnut St.\\
Philadelphia, PA, 19104}
\email{tolya@math.drexel.edu}
\author[Kaliuzhnyi-Verbovetskyi]{Dmitry~S.~Kaliuzhnyi-Verbovetskyi}
\address{Department of Mathematics \\
Drexel University\\
3141 Chestnut St.\\
Philadelphia, PA, 19104}
\email{dmitryk@math.drexel.edu}
\author[Vinnikov]{Victor Vinnikov}
\address{Department of Mathematics\\
Ben-Gurion University of the Negev\\
Beer-Sheva, Israel, 84105} \email{vinnikov@math.bgu.ac.il}
\author[Woerdeman]{Hugo J.~Woerdeman}
\address{Department of Mathematics \\
Drexel University\\
3141 Chestnut St.\\
Philadelphia, PA, 19104}
\email{hugo@math.drexel.edu}

\dedicatory{ Dedicated to Leiba Rodman, a dear friend and wonderful colleague, who
unfortunately passed away too soon}
\thanks{AG, DK-V, HW were partially supported by NSF grant DMS-0901628.
DK-V and VV were partially supported by BSF grant 2010432.}

\subjclass{15A15; 47A13, 13P15, 90C25, 93B28, 47N70}
\keywords{Polynomially defined domain; classical Cartan domains;
contractive realization;
 determinantal representation;
multivariable polynomial; stable polynomial.}

\begin{document}

\begin{abstract} We prove that every matrix-valued rational function $F$, which is regular
on the closure of a bounded domain $\mathcal{D}_\mathbf{P}$ in
$\mathbb{C}^d$ and which has the associated Agler norm strictly
less than 1, admits a finite-dimensional contractive realization
$$F(z)= D + C\mathbf{P}(z)_n(I-A\mathbf{P}(z)_n)^{-1} B. $$
Here $\mathcal{D}_\mathbf{P}$ is
 defined by the inequality $\|\mathbf{P}(z)\|<1$,
where $\mathbf{P}(z)$ is a direct sum of matrix polynomials
$\mathbf{P}_i(z)$ (so that appropriate Archimedean and approximation conditions are
satisfied), and
$\mathbf{P}(z)_n=\bigoplus_{i=1}^k\mathbf{P}_i(z)\otimes I_{n_i}$,
with some
 $k$-tuple $n$ of multiplicities $n_i$; special cases include the open unit polydisk and the classical Cartan domains. The proof uses a matrix-valued version of a Hermitian Positivstellensatz by Putinar,
 and a lurking contraction argument. As a consequence, we show that every polynomial with no zeros on
the closure of $\mathcal{D}_\mathbf{P}$ is a factor of
 $\det (I - K\mathbf{P}(z)_n)$, with a contractive matrix $K$.

\end{abstract}

\maketitle

\section{Introduction}\label{sec:Intro}
It is well-known (see \cite[Proposition 11]{Arov}) that every
 rational matrix function that is contractive on the
open unit disk ${\mathbb D} = \{ z \in {\mathbb C} : |z|<1 \}$ can
be realized as
\begin{equation}\label{real} F(z)= D + zC(I-zA)^{-1} B,
\end{equation} with a contractive (in the spectral norm) colligation matrix $\left[\begin{smallmatrix} A & B \\ C & D
\end{smallmatrix}\right]$. In
several variables, a celebrated result of Agler \cite{Ag} gives
the existence of a realization of the form
\begin{equation}\label{real2} F(z)= D + CZ_\mathcal{X}(I-AZ_{\mathcal{X}})^{-1} B,\qquad
Z_\mathcal{X} = \bigoplus_{i=1}^d z_i I_{{\mathcal X}_i},
 \end{equation}
where $z=(z_1,\ldots,z_d)\in\mathbb{D}^d$ and the colligation
$\left[
\begin{smallmatrix} A & B \\ C & D
\end{smallmatrix}\right]$ is a Hilbert-space unitary operator (with $A$
acting on the orthogonal direct sum of Hilbert spaces ${\mathcal
X}_1, \ldots , {\mathcal X}_d$), for $F$ an operator-valued
function analytic on the unit polydisk $\mathbb{D}^d$ whose Agler
norm $$\|F\|_\mathcal{A}=\sup_{T\in\mathcal{T}}\|F(T)\|\le 1.$$ Here
$\mathcal{T}$ is the set of $d$-tuples $T=(T_1,\ldots,T_d)$ of
commuting strict contractions on a Hilbert space. Such functions
constitute the Schur--Agler class.

Agler's result was generalized to polynomially defined domains in
\cite{AT,BB}. Given a $d$-variable $\ell\times m$ matrix
polynomial $\mathbf{P}$, let
$$\mathcal{D}_\mathbf{P}=\{z\in\mathbb{C}^d\colon\|\mathbf{P}(z)\|<1\},$$
and let ${\mathcal T}_\mathbf{P}$ be the set of $d$-tuples $T$ of
commuting bounded operators on a Hilbert space satisfying
$\|\mathbf{P}(T)\|<1$. Important special cases are:
\begin{enumerate}
\item When $\ell=m=d$ and
    $\mathbf{P}(z)=\diag[z_1,\ldots,z_d]$, the domain
    $\mathcal{D}_\mathbf{P}$ is the unit polydisk
$\mathbb{D}^d$, and  ${\mathcal T}_\mathbf{P}=\mathcal{T}$ is the
set of $d$-tuples of commuting strict contractions.
    \item When $d=\ell m$, $z=(z_{rs})$, $r=1,\ldots,\ell$,
    $s=1,\ldots,m$, $\mathbf{P}(z)=[z_{rs}]$, the domain
    $\mathcal{D}_\mathbf{P}$ is a matrix unit
    ball a.k.a. Cartan's domain of type I. In particular, if $\ell=1$, then
    $\mathcal{D}_\mathbf{P}=\mathbb{B}^d=\{z\in\mathbb{C}^d\colon
    \sum_{i=1}^d|z_i|^2<1\}$ and ${\mathcal T}_\mathbf{P}$ consists
    of commuting strict row contractions $T=(T_1,\ldots,T_d)$.
  \item When $\ell=m$, $d=m(m+1)/2$, $z=(z_{rs})$,  $1\le r\le s\le m$,
 $\mathbf{P}(z)=[z_{rs}]$, where for $r>s$ we set $z_{rs}=z_{sr}$, and  the domain
    $\mathcal{D}_\mathbf{P}$ is a (complex) symmetric matrix unit
    ball a.k.a. Cartan's domain of type II.
      \item When $\ell=m$, $d=m(m-1)/2$, $z=(z_{rs})$,  $1\le r<s\le m$,
 $\mathbf{P}(z)=[z_{rs}]$, where for $r>s$ we set $z_{rs}=-z_{sr}$, and $z_{rr}=0$ for all $r=1,\ldots,m$.
 The domain
    $\mathcal{D}_\mathbf{P}$ is a (complex) skew-symmetric matrix unit
    ball a.k.a. Cartan's domain of type III.
\end{enumerate}

We notice that Cartan domains of 
    types IV--VI can also be represented as $\mathcal{D}_{\mathbf{P}}$, with a linear $\mathbf{P}$; see \cite{Harris} and \cite{Yin}.

For $T \in {\mathcal T}_\mathbf{P}$, the Taylor joint spectrum
$\sigma(T)$ \cite{T1} lies in $\mathcal{D}_\mathbf{P}$ (see
\cite[Lemma 1]{AT}), and therefore for an operator-valued function
$F$ analytic on $\mathcal{D}_\mathbf{P}$ one defines $F(T)$ by
means of Taylor's functional calculus \cite{T2} and
$$ \| F \|_{{\mathcal A},{\mathbf P}}:= \sup_{T
\in {\mathcal T}_{\mathbf{P}}} \| F(T) \|. $$ We say that $F$
belongs to the operator-valued Schur--Agler class associated with
$\mathbf{P}$, denoted by
$\mathcal{SA}_\mathbf{P}(\mathcal{U},\mathcal{Y})$ if $F$ is
analytic on $\mathcal{D}_\mathbf{P}$, takes values in the space
$\mathcal{L}(\mathcal{U},\mathcal{Y})$ of bounded linear operators
from a Hilbert space $\mathcal{U}$ to a Hilbert space
$\mathcal{Y}$, and $\| F \|_{{\mathcal A},{\mathbf P}}\le 1$.

The generalization of Agler's theorem mentioned above that has
appeared first in \cite{AT} for the scalar-valued case and
extended in \cite{BB} to the operator-valued case, says that a
function $F$ belongs to the Schur--Agler class
$\mathcal{SA}_\mathbf{P}(\mathcal{U},\mathcal{Y})$ if and only if
there exists a Hilbert space $\mathcal{X}$ and a unitary
colligation
$$\begin{bmatrix} A & B \\ C & D\end{bmatrix}\colon
(\mathbb{C}^m\otimes \mathcal{X})\oplus
\mathcal{U}\to(\mathbb{C}^\ell\otimes \mathcal{X})\oplus
\mathcal{Y}$$ such that
\begin{equation}\label{eq:realization}
F(z)=D+C(\mathbf{P}(z)\otimes
I_\mathcal{X})\Big(I-A(\mathbf{P}(z)\otimes
I_\mathcal{X})\Big)^{-1}B.
\end{equation}

If the Hilbert spaces $\mathcal{U}$ and $\mathcal{Y}$ are
finite-dimensional, $F$ can be treated as a matrix-valued function
(relative to a pair of orthonormal bases for $\mathcal{U}$ and
$\mathcal{Y}$). It is natural to ask whether every rational
$\alpha\times\beta$ matrix-valued function in the Schur--Agler
class
$\mathcal{SA}_\mathbf{P}(\mathbb{C}^\beta,\mathbb{C}^\alpha)$ has
a realization \eqref{eq:realization} with a contractive
colligation matrix $\left[
\begin{smallmatrix} A & B \\ C & D
\end{smallmatrix}\right]$. This question is open, unless when $d=1$ or
$F$ is an inner (i.e., taking unitary boundary values a.e. on the
unit torus $\mathbb{T}^d=\{z=(z_1,\ldots,z_d)\in\mathbb{C}^d\colon
|z_i|=1,\ i=1,\ldots,d\}$) matrix-valued Schur--Agler function on
$\mathbb{D}^d$. In the latter case, the colligation matrix can be
chosen unitary; see \cite{Knese2011} for the scalar-valued case,
and \cite[Theorem 2.1]{BK} for the matrix-valued generalization.
We notice here that not every inner function is Schur--Agler; see
\cite[Example 5.1]{GKVW} for a counterexample.

 In the present paper, we show that finite-dimensional contractive
realizations of a rational matrix-valued function $F$ exist when
$F$ is regular on the closed domain
$\overline{\mathcal{D}_\mathbf{P}}$ and the Agler norm $\| F
\|_{{\mathcal A},{\mathbf P}}$
 is strictly less than 1 if $\mathbf{P}=\bigoplus_{i=1}^k\mathbf{P}_i$ and the matrix polynomials $\mathbf{P}_i$ satisfy
 a certain natural Archimedean condition. The proof has two
 ingredients: a matrix-valued version of a Hermitian Positivstellensatz \cite{P} (see also \cite[Corollary 4.4]{HP}),
 and a lurking contraction argument. For the first ingredient, we
 introduce the notion of a matrix system of Hermitian quadratic
 modules and the Archimedean property for them, and use the
 hereditary functional calculus for evaluations of a Hermitian
 symmetric matrix polynomial on $d$-tuples of commuting operators
 on a Hilbert space. For the second ingredient, we proceed
 similarly to the lurking isometry argument \cite{Ag,BT,AT,BB},
 except that we are constructing a contractive matrix colligation
 instead of a unitary one.

  We then apply this result to obtain a
determinantal representation $\det(I-K\mathbf{P}_n)$, where $K$ is
a contractive matrix and
$\mathbf{P}_n=\bigoplus_{i=1}^k(\mathbf{P}_i\otimes I_{n_i})$,
with some $k$-tuple $n=(n_1,\ldots,n_k)$ of nonnegative
integers\footnote{We use the convention that if $n_i=0$ then the
corresponding direct summand for $\mathbf{P}_n$ is void.}, for a
multiple of every polynomial which is strongly stable on
${\mathcal{D}_\mathbf{P}}$. (We recall that a polynomial is called
stable with respect to a given domain if it has no zeros in the
domain, and strongly stable if it has no zeros in the domain
closure.) The question of existence of such a representation for a
strongly stable polynomial (without multiplying it with an extra
factor) on a general domain $\mathcal{D}_\mathbf{P}$ is open.

When $\mathcal{D}_\mathbf{P}$ is the open unit polydisk
$\mathbb{D}^d$, the representation takes the form $\det(I-KZ_n)$,
where $Z_n=\bigoplus_{i=1}^dz_iI_{n_i}$,
$n=(n_1,\ldots,n_d)\in\mathbb{Z}_+^d$ (see our earlier work
\cite{GKVW,GKVVW}). In the cases of $\mathbb{D}$ and
$\mathbb{D}^2$, a contractive determinantal representation of a
given stable polynomial always exists; see \cite{Kummert0,GKVVW}. It also exists in the case of 
    multivariable linear functions that are stable on ${\mathbb{D}}^d$, $d=1,2,\ldots$ \cite{GKVW}.
In addition, we showed recently in \cite{GKVVW-PB} that in the
matrix poly-ball case (a direct sum of Cartan domains of type I) a
strongly stable polynomial always has a strictly contractive
realization.

 The paper is organized as follows. In Section
\ref{sec2},  we prove a matrix-valued version of a Hermitian
Positivstellensatz. We then use it in Section \ref{sec:Contr} to
establish the existence of contractive finite-dimensional
realizations for rational  matrix functions from the Schur--Agler
class. In Section \ref{sec:Detrep}, we study contractive
determinantal representations of strongly  stable polynomials.

\section{Positive Matrix Polynomials}\label{sec2}

In this section, we extend the result \cite[Corollary 4.4]{HP} to
 matrix-valued polynomials.
 We will write $A\ge 0$ ($A>0$) when a Hermitian matrix
 (or a self-adjoint operator on a Hilbert space) $A$ is
positive semidefinite (resp., positive definite).
 For a
polynomial with complex matrix coefficients
$$ P(w,z) = \sum_{\lambda,\,\mu} P_{\lambda\mu} w^\lambda z^\mu, $$
where $w=(w_1,\ldots,w_d)$, $z=(z_1,\ldots,z_d)$, and
$w^\lambda=w_1^{\lambda_1}\cdots w_d^{\lambda_d}$, we define
$$ P(T^*,T):= \sum_{\lambda,\,\mu} P_{\lambda\mu} \otimes T^{*\lambda}T^\mu, $$
where $T=(T_1, \ldots , T_d)$ is a $d$-tuple of commuting
operators on a Hilbert space. We will prove that $P$ belongs to a
certain Hermitian quadratic module determined by matrix
polynomials $P_1$, \ldots, $P_k$ in $w$ and $z$ when the
inequalities $P_j(T^*,T)\ge 0$ imply that $P(T^*,T)>0$.

 We denote by $\mathbb{C}[z]$ the
algebra of $d$-variable polynomials with complex coefficients, and
by $\mathbb{C}^{\beta\times \gamma}[z]$ the module over
$\mathbb{C}[z]$ of $d$-variable polynomials with the coefficients
in $\mathbb{C}^{\beta\times \gamma}$. We denote by
$\mathbb{C}^{\gamma\times \gamma}[w,z]_{\rm h}$ the vector space
over $\mathbb{R}$ consisting of polynomials in $w$ and $z$ with
coefficients in $\mathbb{C}^{\gamma\times \gamma}$ satisfying
$P_{\lambda\mu}=P_{\mu\lambda}^*$, i.e., those whose matrix of
coefficients is Hermitian. If we denote by ${P}^*(w,z)$ a
polynomial in $w$ and $z$ with the coefficients $P_{\lambda\mu}$
replaced by their adjoints $P_{\lambda\mu}^*$, then the last
property means that ${P}^*(w,z)=P(z,w)$.

We will say that
$\mathcal{M}=\{\mathcal{M}_\gamma\}_{\gamma\in\mathbb{N}}$ is a
matrix system of Hermitian quadratic modules over $\mathbb{C}[z]$
if the following conditions are satisfied:
\begin{enumerate}
    \item For every $\gamma\in\mathbb{N}$, $\mathcal{M}_\gamma$ is an
    additive subsemigroup of $\mathbb{C}^{\gamma\times \gamma}[w,z]_{\rm
    h}$, i.e., $\mathcal{M}_\gamma+\mathcal{M}_\gamma\subseteq\mathcal{M}_\gamma$.
    \item $1\in\mathcal{M}_1$.
    \item For every $\gamma, \gamma'\in\mathbb{N}$,
    $P\in\mathcal{M}_\gamma$, and $F\in\mathbb{C}^{\gamma\times\gamma' }[z]$, one has
    $F^*(w)P(w,z)F(z)\in\mathcal{M}_{\gamma'}$.
\end{enumerate}
We notice that $\{\mathbb{C}^{\gamma\times \gamma}[w,z]_{\rm
h}\}_{\gamma\in\mathbb{N}}$ is a trivial example of a matrix
system of Hermitian quadratic modules over  $\mathbb{C}[z]$.

\begin{rem}\label{rem:quadmodule}
We first observe that $A\in\mathcal{M}_\gamma$ if
$A\in\mathbb{C}^{\gamma\times \gamma}$ is such that $A=A^*\ge 0$.
Indeed, using (2) and letting $P=1\in\mathcal{M}_1$ and $F$ be a
constant row of size $\gamma$ in (3), we obtain that
$0_{\gamma\times\gamma}\in\mathcal{M}_\gamma$ and that every
constant positive semidefinite $\gamma\times\gamma$  matrix of
rank 1 belongs to $\mathcal{M}_\gamma$, and then use (1). In
particular, we obtain that $I_\gamma\in\mathcal{M}_\gamma$.
Together with (2) and (3) with $\gamma'=\gamma$, this means that
$\mathcal{M}_\gamma$ is a Hermitian quadratic module (see, e.g.,
\cite{PS} for the terminology).

We also observe that, for each $\gamma$, $\mathcal{M}_\gamma$ is a
cone, i.e., it is invariant under addition and multiplication with
positive scalars.

Finally, we observe that $\mathcal{M}$ respects direct sums, i.e.,
$\mathcal{M}_\gamma\oplus\mathcal{M}_{\gamma'}\subseteq\mathcal{M}_{\gamma+\gamma'}$.
In order to see this we first embed $\mathcal{M}_\gamma$ and
$\mathcal{M}_{\gamma'}$ into $\mathcal{M}_{\gamma+\gamma'}$ by
using (3) with $P\in\mathcal{M}_\gamma$, $F=[I_\gamma \
0_{\gamma\times \gamma'}]$ and $P'\in\mathcal{M}_{\gamma'}$,
$F'=[0_{\gamma'\times\gamma}\ I_{\gamma'}]$, and then use (1).
\end{rem}

The following lemma generalizes \cite[Lemma 6.3]{PS}.
\begin{lem}\label{lem:archi}
Let $\mathcal{M}$ be a matrix system of Hermitian quadratic
modules over  $\mathbb{C}[z]$. The following statements are
equivalent:
\begin{itemize}
    \item[(i)] For every $\gamma\in\mathbb{N}$, $I_\gamma$ is an algebraic interior point of
    $\mathcal{M}_\gamma$, i.e.,
    $\mathbb{R}I_\gamma+\mathcal{M}_\gamma=\mathbb{C}^{\gamma\times \gamma}[w,z]_{\rm
h}$.
    \item[(ii)] $1$ is an algebraic interior point of
    $\mathcal{M}_1$, i.e., $\mathbb{R}+\mathcal{M}_1=\mathbb{C}[w,z]_{\rm
h}$.
    \item[(iii)] For every $i=1,\ldots,d$,
     one has $-w_iz_i\in\mathbb{R}+\mathcal{M}_1$.
\end{itemize}
\end{lem}
A matrix system
$\mathcal{M}=\{\mathcal{M}_\gamma\}_{\gamma\in\mathbb{N}}$ of
Hermitian quadratic modules over  $\mathbb{C}[z]$ that satisfies
any (and hence all) of properties (i)--(iii) in Lemma
\ref{lem:archi} is called Archimedean.
{\it Proof.}
(i)$\Rightarrow$(ii) is trivial.

(ii)$\Rightarrow$(iii) is trivial.

(iii)$\Rightarrow$(i). Let
$\mathcal{A}_\gamma=\{F\in\mathbb{C}^{\gamma\times
\gamma}[z]\colon
-F^*(w)F(z)\in\mathbb{R}I_\gamma+\mathcal{M}_\gamma\}$. It
suffices to prove that
$\mathcal{A}_\gamma=\mathbb{C}^{\gamma\times \gamma}[z]$ for all
$\gamma\in\mathbb{N}$. Indeed, any
$P\in\mathbb{C}^{\gamma\times\gamma}[w,z]_{\rm h}$ can be written
as
\begin{multline*}
P(w,z)=\sum_{\lambda,\mu}P_{\lambda\mu}w^\lambda
z^\mu=\row_\lambda[w^\lambda
I_\gamma][P_{\lambda\mu}]\col_\mu[z^\mu I_\gamma] =\\
\row_\lambda[w^\lambda I_\gamma][A^*_\lambda A_\mu-B^*_\lambda
B_\mu]\col_\mu[z^\mu I_\gamma]=A^*(w)A(z)-B^*(w)B(z),
\end{multline*}
where
$$A(z)=\sum_\mu A_\mu z^\mu\in\mathbb{C}^{\gamma\times\gamma}[z],\quad B(z)=\sum_\mu B_\mu
z^\mu\in\mathbb{C}^{\gamma\times\gamma}[z].$$ If
$-B^*(w)B(z)\in\mathbb{R}I_\gamma+\mathcal{M}_\gamma$, then so is
$P(w,z)=A^*(w)A(z)-B^*(w)B(z)$.

By the assumption, $z_i\in\mathcal{A}_{1}$ for all $i=1,\ldots,d$.
We also have  that
$\mathbb{C}^{\gamma\times\gamma}\in\mathcal{A}_\gamma$ for every
$\gamma\in\mathbb{N}$. Indeed, given
$B\in\mathbb{C}^{\gamma\times\gamma}$, we have that
$\|B\|^2I_\gamma-B^*B\ge 0$. By Remark \ref{rem:quadmodule} we
obtain that $\|B\|^2I_\gamma-B^*B\in\mathcal{M}_\gamma$, therefore
$-B^*B\in\mathbb{R}I_\gamma+\mathcal{M}_\gamma$. It follows that
$\mathcal{A}_\gamma=\mathbb{C}^{\gamma\times \gamma}[z]$ for all
$\gamma\in\mathbb{N}$ if $\mathcal{A}_1$ is a ring over
$\mathbb{C}$ and $\mathcal{A}_\gamma$ is a module over
$\mathbb{C}[z]$. We first observe from the identity
\begin{multline*}(F^*(w)+G^*(w))(F(z)+G(z))+(F^*(w)-G^*(w))(F(z)-G(z))=\\ 2(F^*(w)F(z)+G^*(w)G(z))\end{multline*}
for $F,G\in\mathcal{A}_\gamma$ that
\begin{multline*}
-(F^*(w)+G^*(w))(F(z)+G(z))=-2(F^*(w)F(z)+G^*(w)G(z))\\
+(F^*(w)-G^*(w))(F(z)-G(z))\in\mathbb{R}I_\gamma+\mathcal{M}_\gamma,
\end{multline*}
hence $F+G\in\mathcal{A}_\gamma$. Next, for
$F\in\mathcal{A}_\gamma$ and $g\in\mathcal{A}_1$ we can find
positive scalars $a$ and $b$ such that
$aI_\gamma-F^*(w)F(z)\in\mathcal{M}_\gamma$ and
$b-g^*(w)g(z)\in\mathcal{M}_1$. Then we have
\begin{multline*}abI_\gamma-(g^*(w)F^*(w))(g(z)F(z))=\\ b(aI_\gamma-F^*(w)F(z))+F^*(w)\Big((b-g(w)^*g(z))I_\gamma\Big) F(z)\in\mathcal{M}_\gamma,\end{multline*}
Therefore $gF\in\mathcal{A}_\gamma$. Setting $\gamma=1$, we first
conclude that $\mathcal{A}_1$ is a ring over $\mathbb{C}$, thus
$\mathcal{A}_1=\mathbb{C}[z]$. Then, for an arbitrary
$\gamma\in\mathbb{N}$, we conclude that $\mathcal{A}_\gamma$ is a
module over $\mathbb{C}[z]$, thus
$\mathcal{A}_\gamma=\mathbb{C}^{\gamma\times \gamma}[z]$.
\hfill $\square$

Starting with polynomials $P_j\in\mathbb{C}^{\gamma_j\times
\gamma_j}[w,z]_{\rm h}$,
 we introduce the sets ${\mathcal M}_\gamma$, $\gamma \in {\mathbb N}$, consisting of
 polynomials $P\in\mathbb{C}^{\gamma\times\gamma}[w,z]_{\rm h}$
 for which there exist $H_j \in {\mathbb C}^{\gamma_j n_j
\times \gamma} [z]$, for some $n_j \in {\mathbb N}$,  $j=0,\ldots,
k$, such that
\begin{equation}\label{cone}
P(w,z) = H_0^*({w})H_0(z) + \sum_{j=1}^k H_j^*({w}) (P_j(w,z )
\otimes I_{n_j}) H_j(z).
\end{equation}
Here $\gamma_0=1$. We also assume that there exists a constant
$c>0$ such that $c-w_iz_i\in\mathcal{M}_1$ for every
$i=1,\ldots,d$. Then
$\mathcal{M}=\mathcal{M}_{P_1,\ldots,P_k}=\{\mathcal{M}_\gamma\}_{\gamma\in\mathbb{N}}$
is an Archimedean matrix system of Hermitian quadratic modules
generated by $P_1,\ldots,P_k$. It follows from Lemma
\ref{lem:archi} that each $\mathcal{M}_\gamma$ is a convex cone in
the real vector space $\mathbb{C}^{\gamma\times\gamma}[w,z]_{\rm
h}$ and $I_\gamma$ is an interior point in the finite topology
(where a set is open if and only if its intersection with any
finite-dimensional subspace is open; notice that a Hausdorff
topology on a finite-dimensional topological vector space is
unique).

 We can now state the main result of
this section.

\begin{thm} \label{pospol}
Let $P_j\in\mathbb{C}^{\gamma_j\times\gamma_j}[w,z]$, $j=1,\ldots
, k$. Suppose there exists $c>0$ such that $c^2 -w_iz_i\in
{\mathcal M}_1$, for all $i=1,\ldots,d$. Let
$P\in\mathbb{C}^{\gamma\times \gamma}[w,z]$ be such that for every
$d$-tuple $T=(T_1, \ldots , T_d)$ of Hilbert space operators
satisfying $P_j(T^*,T) \ge 0$, $j=1,\ldots , k$, we have that
$P(T^*,T) >0$. Then $P\in {\mathcal M}_\gamma$.
\end{thm}

{\it Proof.}
Suppose that $P \not\in {\mathcal M}_\gamma$. By Lemma
\ref{lem:archi}, $I_\gamma\pm \epsilon P\in\mathcal{M}_\gamma$ for
$\epsilon>0$ small enough. By the Minkowski--Eidelheit--Kakutani
separation theorem (see, e.g., \cite[Section 17]{Kothe}), there
exists a linear functional $L$ on
$\mathbb{C}^{\gamma\times\gamma}[w,z]_{\rm h}$ nonnegative on
${\mathcal M}_\gamma$ such that $L(P) \le 0 < L(I_\gamma)$. For $A
\in {\mathbb C}^{1\times \gamma} [z]$ we define
$$ \langle A, A \rangle = L(A^*(w)A(z)).$$
We extend the definition by polarization:
$$   \langle A, B \rangle = \frac14 \sum_{r=0}^3 i^r \langle A+i^rB,A+i^rB\rangle.
$$
We obtain that $({\mathbb C}^{1\times \gamma}
[z],\langle\cdot,\cdot\rangle)$ is a pre-Hilbert space. Let
${\mathcal H}$ be the Hilbert space completion of the quotient
space
  $ {\mathbb C}^{1\times \gamma} [z] / \{ A \colon \langle A , A \rangle = 0 \}$. Note that $\mathcal H$
  is nontrivial since $L(I_\gamma)>0$.

Next we define multiplication operators $M_{z_i}$, $i=1,\ldots ,
d$, on ${\mathcal H}$. We define $M_{z_i}$ first on the
pre-Hilbert space via $M_{z_i}(A(z)) = z_i A(z)$. Suppose that
$\langle A,A \rangle=0$. Since $c^2 - w_iz_i\in\mathcal{M}_1 $, it
follows that $A^*(w) (c^2 - w_iz_i )A(z) \in {\mathcal M}_\gamma$.
Since $L$ is nonnegative on the cone $\mathcal{M}_\gamma$, we have
\begin{multline*} 0 \le L(A^*(w) (c^2 - w_iz_i)A(z)) =  c^2 \langle A,A \rangle - \langle M_{z_i}(A),M_{z_i}(A)
\rangle =  -\langle M_{z_i}(A),M_{z_i}(A) \rangle.\end{multline*} Thus,
$\langle M_{z_i}(A),M_{z_i}(A) \rangle=0$, yielding that $M_{z_i}$
can be correctly defined on the quotient space. The same
computation as above also shows that $\|M_{z_i}\|\le {c}$ on the
quotient space, and then by continuity this is true on ${\mathcal
H}$. Thus we obtain commuting bounded multiplication operators
$M_{z_i}$, $i=1,\ldots , d$, on ${\mathcal H}$.

Next, let us show that $P_j (M^*,M) \ge 0$, $j=1,\ldots , k$,
where $M=(M_{z_1}, \ldots , M_{z_d} )$. Let $h= [h_r]_{r =
1}^{\gamma_j} \in {\mathbb C}^{\gamma_j} \otimes\mathcal{H}$, and
moreover assume that $h_r$ are elements of the quotient space  $
{\mathbb C}^{1\times \gamma} [z] / \{ A\colon \langle A , A
\rangle = 0 \}$ (which is dense in $\mathcal{H}$).  We will denote
a representative of the coset $h_r$ in
$\mathbb{C}^{1\times\gamma}[z]$ by $h_r(z)$ with a hope that this
will not cause a confusion. Let us compute $\langle P_j(M^*,M) h ,
h \rangle$. We have
$$ P_j (w,z ) = \sum_{\lambda,\,\mu} P_{\lambda\mu}^{(j)}  w^\lambda
z^\mu,\quad P_{\lambda\mu}^{(j)} = [P_{\lambda\mu}^{(j; r,
s)}]_{r, s = 1}^{\gamma_j}.$$ Then $P_j (M^*,M) =
\sum_{\lambda,\,\mu} [P_{\lambda\mu}^{(j; r, s)}M^{* \lambda}
M^\mu]_{r,s = 1}^{\gamma_j}.$ Now \begin{multline*}
 \Big\langle P_j(M^*,M) h, h \Big\rangle = \sum_{r,\,
s =1}^{\gamma_j} \Big\langle \sum_{\lambda,\,\mu}
P_{\lambda\mu}^{(j; r,s)}M^{*\lambda} M^\mu h_r, h_s \Big\rangle
=\\ \sum_{r,\, s =1}^{\gamma_j} \sum_{\lambda,\,\mu}
P_{\lambda\mu}^{(j; r, s)}\Big\langle M^\mu h_r, M^\lambda h_s
\Big\rangle=\\ \sum_{r,\,s =1}^{\gamma_j} \sum_{\lambda,\,\mu}
P_{\lambda\mu}^{(j;r,s)} L \Big( h_s^*(w) w^\lambda z^\mu h_r(z)
\Big)=\\  L \Big( \sum_{r,s=1}^{\gamma_j} h_s^*(w)
\Big(\sum_{\lambda,\,\mu} P_{\lambda\mu}^{(j; r,s)} w^\lambda
z^\mu\Big) h_r(z) \Big) = L(h^*(w) P_j (w,z ) h(z) ),
\end{multline*} which is nonnegative since $h^*(w)P_j (w,z) h(z) \in {\mathcal M}_\gamma$.

By the assumption on $P$ we now have that $P(M^*,M) >0$. By a
calculation similar to the one in the previous paragraph, we
obtain that $L(h^*(w)P(w,z) h(z) ) >0$ for all $h \neq 0$. Choose
now $h(z) \equiv I_\gamma \in \mathbb{C}^{\gamma\times\gamma}[z]$,
and we obtain that $L(P)>0$. This contradicts the choice of $L$.
\hfill $\square$

\section{Finite-dimensional Contractive Realizations}\label{sec:Contr}

In this section, we assume that
$\mathbf{P}(z)=\bigoplus_{i=1}^k\mathbf{P}_i(z)$, where
$\mathbf{P}_i$ are polynomials in $z=(z_1,\ldots,z_d)$ with
$\ell_i\times m_i$ complex matrix coefficients, $i=1,\ldots,k$.
Then, clearly, $\mathcal{D}_\mathbf{P}$ is a cartesian product of
the domains $\mathcal{D}_{\mathbf{P}_i}$. Next, we assume that every $d$-tuple $T$ of commuting
bounded linear operators on a Hilbert space, satisfying $\|  \mathbf{P}(T) \| \le 1$ is a norm limit of elements of ${\mathcal T}_{\mathbf P}$.  We also assume that the
polynomials $P_i(w,z)=I_{m_i}-\mathbf{P}_i^*(w)\mathbf{P}_i(z)$,
$i=1,\ldots,k$, generate an Archimedean matrix system of Hermitian
quadratic modules over $\mathbb{C}[z]$. This in particular means
that the domain $\mathcal{D}_\mathbf{P}$ is bounded, because for
some $c>0$ we have $c^2-w_iz_i\in\mathcal{M}_1$ which implies that
$c^2-|z_i|^2\ge 0$, $i=1,\ldots,d$, when
$z\in\mathcal{D}_\mathbf{P}$. We notice that in the special cases
(1)--(4) in Section \ref{sec:Intro}, the Archimedean condition
holds.

We recall that a polynomial convex hull of a compact set
$K\subseteq\mathbb{C}^d$ is defined as the set of all points
$z\in\mathbb{C}^d$ such that $|p(z)|\le\max_{w\in K}|p(w)|$ for
every polynomial $p\in\mathbb{C}[z]$. A set in $\mathbb{C}^d$ is
called polynomially convex if it agrees with its polynomial convex
hull.
\begin{lem}\label{lem:poly-convex}
$\overline{\mathcal{D}_\mathbf{P}}$ is polynomially convex.
\end{lem}

{\it Proof.}
We first observe that $\overline{\mathcal{D}_\mathbf{P}}$ is
closed and bounded, hence compact. Next, if $z\in\mathbb{C}^d$ is
in the polynomial convex hull of
$\overline{\mathcal{D}_\mathbf{P}}$, then for all unit vectors
$g\in\mathbb{C}^\ell$ and $h\in\mathbb{C}^m$,  one has
$$|g^*\mathbf{P}(z)h|\le\max_{w\in\overline{\mathcal{D}_\mathbf{P}}}|g^*\mathbf{P}(w)h|\le\max_{w\in\overline{\mathcal{D}_\mathbf{P}}}\|\mathbf{P}(w)\|\le
1.$$ Then
$$\|\mathbf{P}(z)\|=\max_{\|g\|=\|h\|=1}|g^*\mathbf{P}(z)h|\le 1,$$
therefore $z\in\overline{\mathcal{D}_\mathbf{P}}$.
\hfill $\square$

\begin{lem}\label{lem:universal-T}
There exists a
$d$-tuple $T_{\rm max}$ of commuting bounded linear operators on a
separable Hilbert space satisfying $\|\mathbf{P}(T_{\rm max})\|\le
1$ and such that
$$\|q(T_{\rm max})\|=\|q\|_{{\mathcal A}, {\bf P}}$$
for every polynomial $q\in\mathbb{C}[z]$.
\end{lem}

{\it Proof.}
The proof is exactly the same as the one suggested in \cite[Page
65 and Exercise 5.6]{Pa} for the special case of commuting
contractions and the Agler norm $\|\cdot\|_\mathcal{A}$ associated
with the unit polydisk, see the first paragraph of Section
\ref{sec:Intro}. Notice that the boundedness of
${\mathcal{D}_\mathbf{P}}$ guarantees that $\|q\|_{{\mathcal A},
{\bf P}}<\infty$ for every polynomial $q$.
\hfill $\square$

\begin{lem}\label{lem:Agler-bdd}
Let $F$ be an
$\alpha\times\beta$ matrix-valued function analytic on
$\overline{\mathcal{D}_\mathbf{P}}$. Then
$\|F\|_{\mathcal{A},\mathbf{P}}<\infty$.
\end{lem}

{\it Proof.}
Since $F$ is analytic on some open neighborhood of the set
$\overline{\mathcal{D}_\mathbf{P}}$ which by Lemma
\ref{lem:poly-convex} is polynomially convex, by the Oka--Weil
theorem (see, e.g., \cite[Theorem 7.3]{AW}) for each scalar-valued
function $F_{ij}$ there exists a sequence of polynomials
$Q^{(n)}_{ij}\in\mathbb{C}^{\alpha\times \beta}[z]$,
$n\in\mathbb{N}$, which converges to $F_{ij}$ uniformly on
$\overline{\mathcal{D}_\mathbf{P}}$. Therefore the sequence of
matrix polynomials $Q^{(n)}=[Q^{(n)}_{ij}]$, $n\in\mathbb{N}$,
converges to $F$ uniformly on $\overline{\mathcal{D}_\mathbf{P}}$.
Let $T$ be any $d$-tuple of commuting bounded linear operators on
a Hilbert space with $\| {\mathbf P} (T) \| \le 1$. By \cite[Lemma 1]{AT}, the
Taylor joint spectrum of $T$ lies in the closed domain
$\overline{\mathcal{D}_{\mathbf{P}}}$ where $F$ is analytic. By
the continuity of Taylor's functional calculus \cite{T2}, we have
that
$$F(T)=\lim_nQ^{(n)}(T).$$
Using Lemma \ref{lem:universal-T}, we obtain that the limit
$$\lim_n\|Q^{(n)}\|_{\mathcal A, {\bf P}}=\lim_n\|Q^{(n)}(T_{\rm max})\|=\|F(T_{\rm max})\|$$
exists and
\begin{multline*}\|F\|_{\mathcal A, {\bf P}}=\sup_{T\in\mathcal{T}_{{\bf P}}}\|F(T)\|
=\sup_{T\in\mathcal{T}_{{\bf P}}}\lim_n\|Q^{(n)}(T)\| \le
\lim_n\|Q^{(n)}(T_{\max})\|=\|F(T_{\rm max})\|<\infty.\end{multline*}
\hfill $\square$

\begin{thm}\label{contr} Let $F$ be a rational $\alpha \times \beta$ matrix function regular on
$\overline{\mathcal{D}_\mathbf{P}}$ and with
 $ \| F
\|_{{\mathcal A},{\mathbf P}} < 1 $. Then there exists
$n=(n_1,\ldots, n_k) \in {\mathbb Z}_+^k$ and a contractive
colligation matrix $\left[
\begin{smallmatrix} A & B \\ C & D
\end{smallmatrix}\right]$ of size $(\sum_{i=1}^k n_im_i + \alpha ) \times
( \sum_{i=1}^k n_i\ell_i + \beta)$ such that
$$ F(z) = D +C\mathbf{P}(z)_n(I-A\mathbf{P}(z)_n)^{-1}B,\qquad \mathbf{P}(z)_n=\bigoplus_{i=1}^k(\mathbf{P}_i(z)\otimes
I_{n_i}).$$
\end{thm}

{\it Proof.} Let $F=QR^{-1}$ with $\det R$ nonzero on
$\overline{\mathcal{D}_\mathbf{P}}$  and let $ \| F \|_{{\mathcal
A},{\mathbf P}} < 1 $. Then we have \begin{equation}\label{eq:est}
R(T)^* R(T)- Q(T)^*Q(T) \ge (1- \| F \|_{{\mathcal A},{\mathbf
P}}^2) R(T)^*R(T) \ge \epsilon^2 I
\end{equation}
 for every
$T\in{\mathcal T}_\mathbf{P}$ with some $\epsilon >0$. Indeed, the
rational matrix function $R^{-1}$ is regular on
$\overline{\mathcal{D}_\mathbf{P}}$. By Lemma \ref{lem:Agler-bdd}
$\|R^{-1}\|_{\mathcal{A},\mathbf{P}}<\infty$. Since
$\|R(T)\|\|R^{-1}(T)\|\ge 1$, we obtain $$\|R(T)
\|\ge\|R^{-1}(T)\|^{-1}\ge\|R^{-1}\|^{-1}_{\mathcal{A},\mathbf{P}}>0,$$
which yields the estimate \eqref{eq:est}.

By Theorem \ref{pospol} there exist $n_0,\ldots,n_k \in{\mathbb
Z}_+$ and polynomials $H_i$ with coefficients in
$\mathbb{C}^{n_im_i\times \beta}$, $i=0,\ldots,k$, (where we set
$m_0=1$) such that by \eqref{cone} we obtain
\begin{multline}\label{AD} R^*(w)R(z)- Q^*(w)Q(z)
 = \\H_0^*(w) H_0(z) +
\sum_{i=1}^kH_i^*(w)
\Big((I-\mathbf{P}_i^*(w)\mathbf{P}_i(z))\otimes I_{n_i}\Big)
H_i(z).\end{multline}
  Denote
\begin{multline*} v(z)= \begin{bmatrix} (\mathbf{P}_1(z)\otimes I_{n_1}) H_1(z) \\
\vdots \\
(\mathbf{P}_k(z)\otimes I_{n_k}) H_k(z)\\
 R(z) \end{bmatrix} \in {\mathbb C}^{(\sum_{i=1}^k \ell_in_i + \beta)\times \beta }[z], x(z)=
 \begin{bmatrix} H_1(z)\\
\vdots\\
H_k(z)\\
 Q(z) \end{bmatrix}  \in {\mathbb C}^{(\sum_{i=1}^k
m_in_i + \alpha)\times \beta }[z]. \end{multline*} Then we may rewrite
\eqref{AD} as
\begin{equation}\label{AD2} v^*(w) v(z) = H_0^*(w) H_0(z) +
x^*(w)x(z). \end{equation} Let us define
$$\mathcal{V}={\rm span} \{ v(z)y\colon z \in {\mathbb C}^d, y\in{\mathbb C}^\beta
\},\quad\mathcal{X}={\rm span} \{ x(z)y\colon z \in\mathbb{C}^d,
y\in{\mathbb C}^\beta\},$$ and let $\{ v(z^{(1)})y^{(1)}, \ldots ,
v(z^{(\nu)})y^{(\nu)} \}$ be a basis for $\mathcal{V} \subseteq
{\mathbb C}^{ \sum_{i=1}^k \ell_i n_i + \beta}$.

{\bf Claim 1.} If $v(z)y = \sum_{i=1}^\nu a_i v(z^{(i)})y^{(i)}$,
then $$x(z)y = \sum_{i=1}^\nu a_i x(z^{(i)})y^{(i)}.$$ Indeed, this
follows from
$$ 0 = \begin{bmatrix} y \\ -a_1y^{(1)} \\ \vdots \\ -a_\nu y^{(\nu)} \end{bmatrix}^*
 \begin{bmatrix} v(z)^* \\ v(z^{(1)})^* \\ \vdots \\
v(z^{(\nu)})^* \end{bmatrix}  \begin{bmatrix} v(z) & v(z^{(1)}) &
\dots & v(z^{(\nu)}) \end{bmatrix}  \begin{bmatrix} y \\
-a_1y^{(1)} \\ \vdots \\ -a_\nu y^{(\nu)} \end{bmatrix} = $$\begin{multline*}
 \begin{bmatrix} y \\ -a_1y^{(1)} \\ \vdots \\ -a_\nu y^{(\nu)} \end{bmatrix}^*
\begin{bmatrix} H_0(z)^*  & x(z)^* \\  H_0(z^{(1)})^*  & x(z^{(1)})^* \\ \vdots & \vdots
\\
H_0(z^{(\nu)})^*  & x(z^{(\nu)})^* \end{bmatrix}  \times \begin{bmatrix}
H_0(z) & H_0(z^{(1)}) & \cdots & H_0(z^{(\nu)})  & \\  x(z) &
x(z^{(1)}) & \dots & x(z^{(\nu)}) \end{bmatrix}  \begin{bmatrix} y
\\ -a_1y^{(1)} \\ \vdots \\ -a_\nu y^{(\nu)} \end{bmatrix} ,\end{multline*} where we
used \eqref{AD2}. This yields
$$  \begin{bmatrix}  H_0(z)  & H_0(z^{(1)}) & \cdots & H_0(z^{(\nu)})
\\
x(z) & x(z^{(1)}) & \dots & x(z^{(\nu)}) \end{bmatrix}
\begin{bmatrix} y \\ -a_1y^{(1)} \\ \vdots \\ -a_\nu y^{(\nu)}
\end{bmatrix} = 0,$$ and thus in particular $x(z)y =
\sum_{i=1}^\nu a_i x(z^{(i)})y^{(i)}$.

We now define  $S\colon \mathcal{V} \to \mathcal{X}$ via
$Sv(z^{(i)})y^{(i)} = x(z^{(i)})y^{(i)}, i=1,\ldots , \nu$. By
Claim 1,
\begin{equation} \label{all} Sv(z)y = x(z)y \ {\rm for\  all}\ z
\in \bigoplus_{j=1}^k {\mathbb C}^{\ell_j\times m_j} \ {\rm and}\
y\in{\mathbb C}^\beta .
\end{equation}

{\bf Claim 2.} $S$ is a contraction. Indeed, let $v =
\sum_{i=1}^\nu a_i v(z^{(i)})y^{(i)} \in \mathcal{V}$. Then $Sv =
\sum_{i=1}^\nu a_i x(z^{(i)})y^{(i)}$, and we compute, using
\eqref{AD2} in the second equality, \begin{multline*}  \| v\|^2 -
\| Sv \|^2 =\\
\begin{bmatrix} a_1y^{(1)} \\ \vdots \\ a_\nu y^{(\nu)}
\end{bmatrix}^* \begin{bmatrix}  v(z^{(1)})^*
 \\ \vdots \\
v(z^{(\nu)})^* \end{bmatrix}  \begin{bmatrix} v(z^{(1)}) & \dots &
v(z^{(\nu)})
\end{bmatrix}  \begin{bmatrix} a_1y^{(1)} \\ \vdots \\ a_\nu
y^{(\nu)}
\end{bmatrix}- \\
\begin{bmatrix} a_1y^{(1)} \\ \vdots \\ a_\nu y^{(\nu)}
\end{bmatrix}^* \begin{bmatrix}  x(z^{(1)})^* \\ \vdots \\
 x(z^{(\nu)})^* \end{bmatrix}  \begin{bmatrix}   x(z) & x(z^{(1)}) & \dots &
x(z^{(\nu)}) \end{bmatrix}  \begin{bmatrix}  a_1y^{(1)} \\ \vdots
\\ a_\nu y^{(\nu)}
\end{bmatrix}=\\ 
 \begin{bmatrix} a_1y^{(1)} \\ \vdots \\ a_\nu y^{(\nu)}\end{bmatrix}^* \begin{bmatrix}  H_0(z^{(1)})^*  \\
 \vdots  \\
 H_0(z^{(\nu)})^*   \end{bmatrix}  \begin{bmatrix}  H_0(z^{(1)}) & \cdots & H_0(z^{(\nu)}) \end{bmatrix}
 \begin{bmatrix} a_1y^{(1)}
 \\ \vdots \\ a_\nu y^{(\nu)} \end{bmatrix} \ge 0,
 \end{multline*}
proving Claim 2.

Extending $S$ to the contraction
$S_{ext}=\left[\begin{smallmatrix} A & B \\ C & D
\end{smallmatrix}\right]\colon {\mathbb C}^{ \sum_{i=1}^k\ell_i
n_i + \beta} \to {\mathbb C}^{ \sum_{i=1}^k m_in_i + \alpha}$ by
setting $S_{ext} |_{\mathcal{V}^\perp} = 0$, we
 now obtain from \eqref{all} that
$$ A\mathbf{P}(z)_nH(z) + B R(z) = H(z) , \ \ C \mathbf{P}(z)_nH(z) + DR(z)=
Q(z). $$ Eliminating $H(z)$, we arrive at
$$ (D+ C\mathbf{P}(z)_n (I- A\mathbf{P}(z)_n)^{-1}B)R(z) = Q(z), $$
yielding the desired realization for $F=QR^{-1}$.
\hfill $\square$

The following statement is a special case of Theorem \ref{contr}.
\medskip

\begin{cor} Let $F$ be a rational matrix function regular on the closed bidisk
$\overline{\mathbb{D}^2}$ such that
$$\|F\|_\infty=\sup_{(z_1,z_2) \in {\mathbb D}^2}
 \|F (z_1,z_2) \| <1.$$
Then $ F$ has a finite dimensional contractive realization
\eqref{real2}, that is, there exist $n_1,n_2\in\mathbb{Z}_+$ such
that $\mathcal{X}_i=\mathbb{C}^{n_i}$, $i=1,2$, and
$Z_\mathcal{X}=z_1I_{n_1}\oplus z_2I_{n_2}=Z_n$.
\end{cor}

{\it Proof.} One can apply Theorem \ref{contr} after observing that on the bidisk the Agler norm and the supremum norm
coincide, a result that goes back to \cite{Ando}. \hfill $\square$

\section{Contractive Determinantal Representations}\label{sec:Detrep}
Let a polynomial $\mathbf{P}=\bigoplus_{i=1}^k\mathbf{P}_i$  and a
domain $\mathcal{D}_\mathbf{P}$ be as in Section \ref{sec:Contr}.
We apply Theorem \ref{contr} to obtain a contractive determinantal
representation for a multiple of every polynomial strongly stable
on $\mathcal{D}_\mathbf{P}$. Please notice the analogy with the main result in \cite{Kummer}, where a similar result is obtained in the setting of definite determinantal representation for hyperbolic polynomials

\begin{thm}\label{thm:sssr}
Let $p$ be a polynomial in $d$ variables $z=(z_1,\ldots,z_d)$,
which is strongly stable on $\mathcal{D}_\mathbf{P}$. Then there
exists a polynomial $q$, nonnegative integers $n_1,\ldots,n_k$,
and a contractive
 matrix
 $K$ of size $\sum_{i=1}^km_in_i\times\sum_{i=1}^k\ell_in_i$ such that
$$ p(z)q(z) = \det (I - K\mathbf{P}(z)_n),\qquad  \mathbf{P}(z)_n=\bigoplus_{i=1}^k(\mathbf{P}_i(z)\otimes I_{n_i}).$$
\end{thm}

{\it Proof.}
Since $p$ has no zeros in $\overline{\mathcal{D}_{\bf P}}$, the rational function
$g=1/p$ is regular on $\overline{\mathcal{D}_{{\bf P}}}$. By
Lemma \ref{lem:Agler-bdd}, $\| g \|_{{\mathcal A}, {\bf P}} <
\infty$. Thus we
 can find a constant $c>0$ so that  $\|cg\|_{{\mathcal A}, {\bf P}} <1$.
Applying now Theorem \ref{contr}
 to $F=cg$, we obtain a $k$-tuple $n=(n_1,\ldots,n_k)\in\mathbb{Z}_+^k$ and a
 contractive colligation matrix $\left[\begin{smallmatrix} A & B \\ C & D \end{smallmatrix}\right]$ so that

\begin{multline}\label{contr-real}
cg(z)=\frac{c}{p(z)}
 = D+C{\bf P}(z)_n(I-A{\bf P}(z)_n)^{-1}B
=
 \frac{\det \begin{bmatrix}  I-A{\bf P}(z)_n & B \\ -C{\bf P}(z)_n &
 D \end{bmatrix}}{\det (I-A{\bf P}(z)_n )}.
\end{multline}
This shows that $$\frac{\det (I-A{\bf P}(z)_n )}{p(z)}$$
is a polynomial. Let  $K= A$. Then $K$ is a 
contraction, and $$q(z)=\frac{\det (I-K{\bf P}(z)_n )}{p(z)}$$ is
a polynomial.
\hfill $\square$

\begin{rem}\label{rem:q-ss}
Since the polynomial $\det (I-K{\bf P}(z)_n)$
in Theorem \ref{thm:sssr} is stable on
$\mathcal{D}_\mathbf{P}$, so is $q$.
\end{rem}


\begin{thebibliography}{10}

\bibitem{Ag} J.~Agler.
\newblock On the representation of certain holomorphic functions
defined on a polydisc,
\newblock In  Topics in operator theory: Ernst D. Hellinger Memorial
Volume,  {\em Oper. Theory Adv. Appl.},
 Vol.~\textbf{48}, pp. 47--66,
  Birkh\"auser, Basel, 1990.

\bibitem{AW}
 H. Alexander and J. Wermer. {\em Several complex variables and Banach algebras.} Third edition.
 Graduate Texts in Mathematics, 35. Springer-Verlag, New York, 1998.


\bibitem{AT}
C.-G. Ambrozie and D. A. Timotin. Von Neumann type inequality for
certain domains in
 $\mathbb{C}^n$. {\em Proc. Amer. Math. Soc.} 131 (2003), no. 3, 859--869 (electronic).

\bibitem{Ando} T. Ando. On a pair of commutative contractions. {\em Acta Sci. Math. (Szeged)} 24 (1963), 88--90.

\bibitem{Arov} D. Z. Arov.
Passive linear steady-state dynamical systems. (Russian) {\em
Sibirsk. Mat. Zh.} 20 (1979), no. 2, 211--228, 457.


 \bibitem{BB}
 J. A. Ball and V. Bolotnikov. Realization and interpolation for Schur--Agler-class functions on domains with matrix
 polynomial defining function in $\mathbb{C}^n$. {\em J. Funct. Anal.} 213 (2004), no. 1, 45--87.

\bibitem{BK}  J. A. Ball and D. S. Kaliuzhnyi-Verbovetskyi.
  Rational Cayley inner Herglotz--Agler functions: Positive-kernel decompositions and
   transfer-function realizations. {\em Linear Algebra Appl.}  456 (2014), 138--156.

\bibitem{BT} J. A. Ball and T. T. Trent. Unitary colligations,
reproducing kernel Hilbert spaces, and Nevanlinna--Pick
interpolation in several variables. {\em J. Funct. Anal.} 157
(1998), no. 1, 1--61.

\bibitem{GKVW} A. Grinshpan, D. S. Kaliuzhnyi-Verbovetskyi, and H. J. Woerdeman.
Norm-constrained determinantal representations of multivariable
polynomials. {\em Complex Anal. Oper. Theory} 7 (2013), 635--654.


\bibitem{GKVVW} A. Grinshpan, D. S. Kaliuzhnyi-Verbovetskyi, V. Vinnikov, and H. J. Woerdeman.
Stable and real-zero polynomials in two variables. {\em Multidim.
Syst. Sign. Process.} 27 (2016), 1-26.

\bibitem{GKVVW-PB} A. Grinshpan, D. S. Kaliuzhnyi-Verbovetskyi, V. Vinnikov, and H. J. Woerdeman. Contractive determinantal representations of stable polynomials
on a matrix polyball.  {\it Math. Z.} 283 (2016), 25--37

\bibitem{Harris} L. A. Harris. Bounded symmetric homogeneous domains in infinite dimensional spaces. In ``Proceedings on Infinite Dimensional Holomorphy" (Internat. Conf., Univ. Kentucky, Lexington, Ky., 1973), pp. 13--40. {\it Lecture Notes in Math.}, Vol. 364, Springer, Berlin, 1974.

\bibitem{HP} J. W. Helton and M. Putinar. Positive polynomials in scalar and matrix variables, the spectral theorem, and optimization.
In {\em Operator theory, structured matrices, and dilations}, 229--306, Theta Ser. Adv. Math., 7, Theta, Bucharest, 2007.


\bibitem{Knese2011} G. Knese. Rational inner functions in the Schur-Agler class of the polydisk. {\em Publ. Mat.},
55 (2011), 343--357.


\bibitem{Kothe} G. K\"othe. Topologische lineare R\"aume. I. (German) Zweite verbesserte Auflage.
 Die Grundlehren der Mathematischen Wissenschaften, Band 107 Springer-Verlag, Berlin-New York 1966

\bibitem{Kummer} M. Kummer, Determinantal Representations and the B\'ezout Matrix, arXiv:1308.5560.

\bibitem{Kummert0}
A. Kummert.
\newblock Synthesis of two-dimmensional lossless
$m$-ports with prescribed scattering matrix.
\newblock {\em Circuits Systems Signal Processing} 8 (1989), no. 1, 97--119.


\bibitem{P}
M. Putinar. On Hermitian polynomial optimization. {\em Arch. Math.
(Basel)} 87 (2006), no. 1, 41--51.

\bibitem{Pa}
V.~Paulsen. {\em Completely bounded maps and operator algebras.}
Cambridge Studies in Advanced Mathematics, 78. Cambridge
University Press, Cambridge, 2002.

\bibitem{PS}
M. Putinar and C. Scheiderer. Quillen property of real algebraic
varieties. {\it M\"unster J. of Math.} 7 (2014), 671--696


\bibitem{T1}
J. L. Taylor.  A joint spectrum for several commuting operators.
{\em J. Functional Analysis} 6 (1970), 172--191.

\bibitem{T2}
J. L. Taylor. The analytic-functional calculus for several
commuting operators. {\em Acta Math.} 125 (1970), 1--38.

\bibitem{Yin} W. Yin. Two problems on Cartan domains.  Preprint, 2006,  arXiv:0603205.

\end{thebibliography}
\end{document}